\documentclass[a4paper,11pt]{article}

\usepackage[utf8]{inputenc} 
\usepackage[T1]{fontenc} 
\usepackage{amsmath}
\usepackage{amsfonts}
\usepackage{amsthm}
\usepackage{amssymb}
\usepackage{color}
\usepackage{graphicx}
\usepackage{float} 
\usepackage{subfigure}
\usepackage{epsfig}
\usepackage{url}

\usepackage[superscript]{cite} 
\usepackage{setspace} 
\makeatletter
\renewcommand\@biblabel[1]{#1.}
\makeatother

\def\ds{\displaystyle }

\def\x{{\bf x}}
\def\u{{\bf u}}
\def\d{{\bf d}}
\def\B{{\bf B}}

\def\Gp{\Gamma_p }
\def\Go{\Gamma_O }
\def\Gi{\Gamma_I }

\graphicspath{{./fig/}{./fig/west/}{./fig/aug/}{./fig/jet/}} 
\makeatletter
\def\input@path{{.}{./fig/}}
\makeatother

\begin{document}

\begin{center}
	{\Large TOKAMAK PLASMA BOUNDARY RECONSTRUCTION USING TOROIDAL HARMONICS AND AN OPTIMAL CONTROL METHOD }\\
	{Blaise Faugeras}\\
{\footnotesize CASTOR Team-Project, INRIA and Laboratoire J.A. Dieudonn\'e CNRS UMR 7351}\\ 
{\footnotesize Universit\'e Nice Sophia-Antipolis, Parc Valrose, 06108 Nice Cedex 02, France}\\
{\footnotesize E-mail: Blaise.Faugeras@unice.fr}
\end{center}


%


\section*{Abstract}
This paper proposes a new fast and stable algorithm for the reconstruction of the plasma boundary 
from discrete magnetic measurements taken at several locations 
surrounding the vacuum vessel. The resolution of this inverse problem takes two steps. 
In the first one we transform the set of measurements 
into Cauchy conditions on a fixed contour $\Gamma_O$ close to the measurement points. 
This is done by least square fitting a truncated series of toroidal harmonic functions to the 
measurements. The second step consists in solving a Cauchy problem for the elliptic equation 
satisfied by the flux in the vacuum and for the overdetermined boundary conditions on $\Gamma_O$ 
previously obtained with the help of toroidal harmonics. 
It is reformulated as an optimal control problem on a fixed annular domain of external boundary $\Gamma_O$ 
and fictitious inner boundary $\Gamma_I$. A regularized Kohn-Vogelius 
cost function depending on the value of the flux on $\Gamma_I$ and measuring the discrepency 
between the solution to the equation satisfied by the flux obtained using Dirichlet conditions on $\Gamma_O$ 
and the one obtained using Neumann conditions is minimized. 
The method presented here has led to the development of a software, called VacTH-KV, 
which enables plasma boundary reconstruction in any Tokamak. \\
{\bf Keywords}: inverse problem, boundary identification, magnetic measurements


\section*{\bf I. INTRODUCTION}
In order to control the plasma during a fusion experiment in a Tokamak 
it is important to be able to compute its boundary in the vacuum vessel. 
This boundary is deduced from the knowledge of the poloidal flux which itself  
is computed from a set of discrete magnetic measurements of the poloidal magnetic field and flux 
scattered around the vacuum vessel.

This paper presents a new fast and stable algorithm for the reconstruction of the poloidal flux in the vacuum surrounding the 
plasma and of the plasma boundary.

Let us first briefly recall the equations modelizing the equilibrium of a plasma in a Tokamak \cite{Wesson:2004}.

Assuming an axisymmetric configuration one considers a 2D poloidal cross section of the vacuum vessel 
in the $(r,z)$ system of coordinates. Throughout the manuscript bold face letters represent vectors or matrices. 
In this setting the poloidal flux $\psi(r,z)$ is linked to the magnetic field through 
the relation $\mathbf{B}=(B_r,B_z)=\ds\frac{1}{r}(-\frac{\partial \psi}{\partial z},\frac{\partial \psi}{\partial r})$ 
and, as there is no toroidal current density in the vacuum outside the plasma, satisfies the following linear 
elliptic partial differential equation   
\begin{equation*}
\label{eqn:psi}
L \psi=0 
\end{equation*}
 where $L$ denotes the elliptic operator 
$$
L.= -[\ds \frac{\partial }{\partial r}(\frac{1}{r}\frac{\partial .}{\partial r}) 
+ \ds \frac{\partial }{\partial z}(\frac{1}{r}\frac{\partial .}{\partial z}) ]
$$ 

\begin{figure}[h]
\centering
\def\svgwidth{0.7\linewidth}
\input{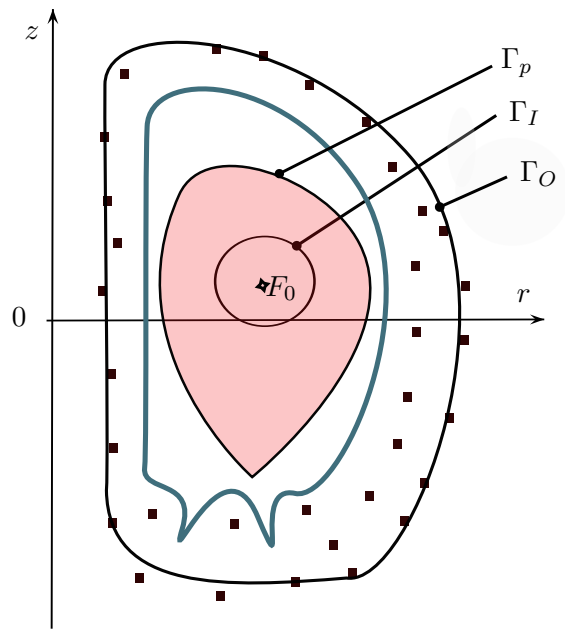}
\caption{Skematic representation of the poloidal cross section of a Tokamak. 
	The unknown plasma domain of boundary $\Gamma_p$ and the limiter contour are shown. 
	Magnetic sensors are depicted by squares. The external contour $\Gamma_O$, 
	the inner contour $\Gamma_I$ as well as the current center $F_0$ are represented.}
\label{fig:geom}
\end{figure}

The unknown plasma boundary $\Gamma_p$ is determined from the equation \\
$\psi(r,z)=\psi_p$, 
$\psi_p$ being the value of the flux at an X-point or the value of the flux for the outermost 
closed flux line inside the limiter (see Figure \ref{fig:geom}). 
The reconstruction of this unknown plasma boundary is a Cauchy problem for 
the operator $L$. A supplementary difficulty is that one does not have directly access 
to a complete set of Cauchy data on an external measurement countour $\Gamma_O$ but only 
to discrete magnetic measurements from sensors not necessarily positioned 
on a single contour.

The proposed resolution method for this inverse problem takes two steps. 
In the first one we transform the set of discrete measurements 
into Cauchy conditions, $\psi$ and its normal derivative, on a fixed contour $\Gamma_O$ 
chosen close to the measurement points. This is done by least square fitting a truncated series 
of toroidal harmonic functions to the magnetic measurements. In principle this first step can be sufficient to 
determine the plasma boundary 
\cite{Lee:1981,Deshko:1983,Bondarenko:1984,Alladio:1986, Kurihara:1992, VanMilligen:1994,
Fischer:2011,ACL.B.Faugeras.14.2} (also see the review article \cite{Braams:1991}). 
Indeed the toroidal harmonics expansion is valid anywhere in the 
vacuum surrounding the plasma. The method of toroidal harmonics presented in \cite{ACL.B.Faugeras.14.2} 
proved to be very efficient for the WEST Tokamak configuration \cite{Bucalossi:2011}. 
The equivalent of magnetic measurements were generated with the 
equilibrium code CEDRES++ \cite{ACL.B.Faugeras.15.1}. 
Even when noise is added to the synthetic measurements the boundary reconstructions are accurate 
with a judicious choice of the order of the toroidal harmonics used.

However since then new numerical experiments conducted for other Tokamaks with more elongated plasma shapes than WEST and with real 
experimental measurements showed that due to the ill-posedness of the inverse problem, the plasma boundary reconstructed 
after this first step involving toroidal harmonics only can in some cases be inaccurate. 
Nevertheless the interpolation of the discrete magnetic measurements on a contour $\Gamma_O$ chosen to be close to 
the measurement points is allways very accurate even though the extrapolation of the flux 
towards the plasma boundary can suffer from some perturbations due to the singularity 
of internal toroidal harmonics. This is the reason why we propose the following second step which can 
be seen as a regularization procedure for the toroidal harmonics expansion.

The second step consists in solving a Cauchy problem for the elliptic equation $L \psi =0$ 
satisfied by the flux and for the overdetermined boundary conditions on $\Gamma_O$ 
obtained with the help of toroidal hamonics as explained above. 
The proposed method consists in a reformulation of the problem as an optimal control problem on 
a fixed annular domain of external boundary $\Gamma_O$ and an fictitious inner boundary $\Gamma_I$ (see Figure \ref{fig:geom}) 
located inside the plasma domain. 
We aim at minimizing a regularized Kohn-Vogelius \cite{Kohn:1985} cost function depending on the value of the flux 
on $\Gamma_I$ and measuring the discrepency between the solution to the 
equation satisfied by $\psi$ obtained using Dirichlet conditions on $\Gamma_O$ 
and the one obtained using Neumann conditions. 

The analysis of this second step alone was presented 
by the author in \cite{ACL.B.Faugeras.12.2} but only tested on a very academic case with synthetic 
Cauchy data. The algorithm proposed in this paper is a combination of the two previous works 
\cite{ACL.B.Faugeras.14.2} and \cite{ACL.B.Faugeras.12.2}. 
It combines the precise data interpolation provided 
by the toroidal harmonics method to generate Cauchy conditions on $\Gamma_O$ and 
the robustness of the optimal control approach to extrapolate these Cauchy conditions towards the plasma 
and compute the plasma boundary. An implementation for fast numerical computations is also proposed.
 
The introduction of an optimal control problem on fixed domain 
appears in \cite{Blum:1987, Blum:1989} where the cost function however is not of the same 
type as the one proposed in this paper. 
The advantage of the cost function proposed in this work is that it uses the 
given Dirichlet and Neumann boundary conditions on the external contour $\Gamma_O$ in a symmetric way. 
For the analysis of section II.B to be valid the control variable has to be the Dirichlet conditions 
on the internal contour $\Gamma_I$.

The use of a fixed inner contour also appears in \cite{Feneberg:1984} 
where a surface current sheet modelizes the plasma and in \cite{Hakkarainen:1984, Kurihara:1993, Kurihara:2000} 
in the framework of boundary integral equations which might not be easily applicable 
in the case of iron core Tokamaks.

The next section presents the two main steps of the proposed method implemented in the code VacTH-KV. 
Section III
then gives details on the numerical methods used and shows some numerical results.

\section*{\bf II. OVERVIEW OF THE METHOD}

\subsection*{\bf II.A. Step1: interpolation of discrete magnetic measurements with toroidal harmonics}

The poloidal flux at any point $\mathbf{x}=(r,z)$ of the vacuum surrounding the plasma is represented by
\begin{equation*}
	\psi(\x)=\psi_C(\x) + \psi_{TH}(\x)
\end{equation*}
The first term $\psi_C$  represents the flux generated by the poloidal field coils. It is computed 
using Green functions. Each coil is modelized by a number of filaments of current of given intensity.

The second term $\psi_{TH}$ as a solution to the equation $L \psi_{TH} = 0$ 
can be uniquely decomposed in a series of toroidal harmonics. 
These toroidal harmonics are found from computations \cite{Braams:1986,Fischer:2011} 
involving a quasi separation of variable technique which lead to define external and internal harmonics.
They form a complete set of solutions to this equation in any annular domain \cite{Fischer:2011}. 
In particular this implies that this expansion can be used for an iron core Tokamak such as WEST.

The $\psi_{TH}$ term is approximated by the following truncated expansion 
\begin{equation*}
\left\{
\begin{array}{rcl}
\psi_{TH}&=&\psi_{ext} + \psi_{int}\\
\psi_{ext}&=& \ds  \frac{r_0 \sinh \zeta}{\sqrt{\cosh \zeta - \cos \eta}} [ \ds \sum_{n=0}^{n_e}(a_n^e  \cos(n\eta) + b_n^e \sin(n\eta))Q^1_{n-1/2}(\cosh\zeta)  ] \\
\psi_{int}&=& \ds \frac{r_0 \sinh \zeta}{\sqrt{\cosh \zeta - \cos \eta}} [ \ds \sum_{n=0}^{n_i} (a_n^i  \cos(n\eta) +  b_n^i  \sin(n\eta))P^1_{n-1/2}(\cosh\zeta) ] \\
\end{array}
\right.
\end{equation*}
The toroidal harmonic functions are computed in a toroidal coordinate system $(\zeta,\eta)$ 
which depends on the choice of a pole $F_0$. 
The internal harmonics involve the Legendre functions 
of the first kind, of degree one and half-integer order, $P^1_{n-1/2}$ which are singular at the pole $F_0$.  
The external harmonics involve the Legendre functions 
of the second kind, of degree one and half-integer order, $Q^1_{n-1/2}$ which are singular on the axis $r=0$.
If we denote by  
$\u=(a_0^e, \dots , a_{n^e}^e, b_1^e, \dots , b_{n^e}^e,a_0^i, \dots , a_{n^i}^i, b_1^i, \dots , b_{n^i}^i)$
the vector of the coefficients of the expansion we have an analytic expression for $\psi_{TH}(\x;\u)$. 
Its gradient can be computed analytically as well leading to an expression for 
the magnetic field in terms of $\u$, $\mathbf{B}_{TH}(\x;\u)$. It remains to identify $\u$ from the set of magnetic 
measurements available. These discrete magnetic measurements are of three types:
\begin{itemize}
	\item B probes provide $N_B$ measurements of the poloidal field at points $\x_i^B$ 
		and directions given by unit vectors $\d_i$ such that $B_i^{meas} \approx \B(\x_i^B).\d_i$ 
		where the dot represents the scalar product
	\item Flux loops provide $N_f$ flux measurements at points $\x_i^f$ such that $\psi_i^{meas} \approx \psi(\x_i^f)$  
	\item Saddle loops provide $N_s$ flux variations measurements between two points such that 
		$\delta_i\psi^{meas} \approx \psi(\x_i^1)-\psi(\x_i^2) := \delta \psi(\x_i^1,\x_i^2) $ 
\end{itemize}

The unknown vector $\u$ of coefficients is computed minimizing the least square cost function 
\begin{equation*}
\left\{
\begin{array}{lcl}
	J(\u)&=&\ds \sum_{i=1}^{N_B}\frac{(B_{TH}(\x_i^B;\u) - \tilde{B}^{meas}_i)^2}{\sigma_B^2}\\
	&&+\ds \sum_{i=1}^{N_f}\frac{(\psi_{TH}(\x_i^f;\u) - \tilde{\psi}^{meas}_i)^2}{\sigma_f^2} 
	+\ds \sum_{i=1}^{N_s}\frac{( \delta \psi_{TH}(\x_i^1,\x_i^2;\u) - \delta_i\tilde{\psi}^{meas})^2}{\sigma_s^2}  
\end{array}
\right.
\label{eqn:costfunction}
\end{equation*}
Where $B_{TH}(\x_i^B;\u)=\B_{TH}(\x_i^B;\u).\d_i$ 
and the known contributions of the poloidal field coils are substracted from the measurements:
\begin{equation*}
\left\{
\begin{array}{lcl}
\tilde{B}_i^{meas}&=&B_i^{meas} - {\B}_C(\x_i^B).\d_i,\quad \mathrm{for} \  i=1, \dots N_B\\
\tilde{\psi}_i^{meas}&=& \psi_i^{meas} - {\psi}_C(\x_i^f),\quad \mathrm{for} \ i=1, \dots N_f\\
\delta_i\tilde \psi^{meas}&=&\delta_i \psi^{meas} - \delta \psi_C(\x_i^1, \x_i^2),\quad \mathrm{for} \ i=1, \dots N_s
\end{array}
\right.
\end{equation*}
The weights $\sigma_B$ and $\sigma_f$, $\sigma_s$ correspond to the assumed measurement errors. 
$J(\u)$ is quadratic and is minimized by solving the normal equation giving the optimal set of coefficients $\u_{opt}$.

Once $\u_{opt}$ is computed an approximation of the flux 
can be obtained at any point of the vacuum surrounding the plasma 
by 
$$
\psi(\x;\u_{opt})=\psi_{TH}(\x;\u_{opt})+\psi_C(\x)
$$

In particular one can evaluate $\psi$ and its normal derivative 
in order to provide Cauchy boundary conditions on any fixed closed countour $\Gamma_O$. 
In this work $\Gamma_O$ is a contour chosen to be close to the measurements. If all 
sensors were located a single contour this sensor or measurement contour 
would be a natural choice for $\Gamma_O$.

\subsection*{\bf II.B. Step 2: an optimal control method}
\label{subsec:step2}

From the method described in the previous paragraph one can compute a complete set of 
Cauchy conditions, $f=\psi$ on $\Gamma_O$ 
and $g = \ds \frac{1}{r}\frac{\partial \psi}{\partial n}$ on $\Gamma_O$. 
It is then possible to employ the optimal 
control method discussed in \cite{ACL.B.Faugeras.12.2} which is summed up below.

The poloidal flux satisfies 
\begin{equation*}
\label{eqn:cauchypbm0}
\left \lbrace
\begin{array}{l}
L \psi = 0 \quad \mathrm{in}\ \Omega_X \\[10pt]
\psi= f \quad  \mathrm{on}\  \Gamma_O \\[10pt]
\ds \frac{1}{r}\frac{\partial \psi}{\partial n}= g\quad  \mathrm{on}\  \Gamma_O \\[10pt]
\psi=\psi_p \quad \mathrm{on}\ \Gamma_p 
\end{array} 
\right.
\end{equation*}

In this formulation the annular domain $\Omega_X=\Omega_X(\psi)$ 
contained between $\Gamma_O$ and $\Gamma_p$ (see Fig. \ref{fig:geom}) is unknown since 
the free plasma boundary $\Gamma_p$ is unknown. 
Moreover the problem is ill-posed as there are two Cauchy conditions on the boundary $\Gamma_O$.
 
In order to compute the flux in the vacuum and to find the plasma boundary we define a new problem 
approximating the original one. 
Let us define a fictitious boundary $\Gamma_I$ fixed inside the plasma. 
We seek an approximation of the poloidal flux 
$\psi$ satisfying $L\psi=0$ in the domain $\Omega$ contained between the fixed boundaries 
$\Gamma_O$ and $\Gamma_I$. The problem becomes one formulated on the fixed annular domain $\Omega$.
Freezing the domain to $\Omega$ by introducing the fictitious boundary $\Gamma_I$ 
enables to remove the nonlinearity of the problem. The plasma boundary $\Gamma_p$ 
can still be computed as an iso-flux line and thus is an output of our computations.

Let us insist here on the fact that this problem is an approximation to the original one since in 
the domain between $\Gamma_p$ and $\Gamma_I$, $\psi$ should satisfy the Grad-Shafranov equation. 
The relevance of this approximating model is consolidated by the Cauchy-Kowalewska theorem \cite{Courant:1962}. 
For $\Gamma_p$ smooth enough the function $\psi$ can be extended in the sense 
of $L\psi=0$ in a neighborhood of $\Gamma_p$ inside the plasma. 
Hence the problem formulated on a fixed domain with a fictitious boundary 
$\Gamma_I$ not "too far" from $\Gamma_p$ is an approximation of the free boundary problem. 
If $\Gamma_I$ were identical with $\Gamma_P$ 
then by the virtual shell principle \cite{Shafranov:1972} 
the quantity $w=\ds \frac{1}{r}\frac{\partial \psi}{\partial n}|_{\Gamma_I}$ would 
represent the surface current density (up to a factor $\ds \frac{1}{\mu_0}$) 
on $\Gamma_p$ for which the magnetic field created outside the plasma by the current 
sheet is identical to the field created by the real current density spread throughout the plasma.

The shape of the internal contour $\Gamma_I$ is chosen to be a circle or an ellipse. 
Although no thorough sensitivity analysis has been conducted concerning the shape and size of $\Gamma_I$ 
the numerical results presented in the next section did not prove to be very dependent on these factors.

No boundary condition is known on $\Gamma_I$ and boundary conditions are over determined on $\Gamma_O$. 
One way to deal with this and to solve such a problem 
is to formulate it as an optimal control one. 
We are going to compute a function $\psi$ such that 
the Dirichlet boundary condition $u=\psi$ on $\Gamma_I$ is such 
that the Cauchy conditions on $\Gamma_O$ are satisfied 
as nearly as possible in the sense of the Kohn-Vogelius error functional defined below which 
enables to use Dirichlet and Neumann conditions in a symmetric way.

We consider two separate well-posed sub-problems. In the first one we retain the Dirichlet boundary condition on 
$\Gamma_O$ only, assume $v$ is given on $\Gamma_I$ and seek the solution $\psi_D$ 
of the boundary value problem:

\begin{equation*}
\label{eqn:psid}
\left\{
\begin{array}{l}
L \psi_D = 0 \quad \mathrm{in}\ \Omega \\[10pt]
\psi_D = f \quad  \mathrm{on}\  \Gamma_O \\[10pt]
\psi_D = v\quad  \mathrm{on}\  \Gamma_I 
\end{array}
\right.
\end{equation*}

This solution $\psi_D$ 
can be decomposed  in a part linearly depending on $v$ 
and a part depending on $f$ only. 
We have the following decomposition:
\begin{equation*}
\label{eqn:decompositionD}
\psi_D=\psi_D(v,f)=\psi_D(v,0)+\psi_D(0,f):=\psi_D(v)+\tilde{\psi}_D(f)
\end{equation*}

In the second sub-problem we retain the Neumann boundary condition only 
and look for $\psi_N$ satisfying the boundary value problem:
\begin{equation*}
\label{eqn:psin}
\left\{
\begin{array}{l}
L \psi_N = 0 \quad \mathrm{in}\ \Omega \\[10pt]
\ds \frac{1}{r}\frac{\partial \psi_N}{\partial n}= g\quad  \mathrm{on}\  \Gamma_O \\[10pt]
\psi_N= v \quad  \mathrm{on}\  \Gamma_I 
\end{array}
\right.
\end{equation*}

in which $\psi_N$  can be decomposed  in a part linearly depending on $v$ 
and a part depending on $g$ only. 
We have the following decomposition:
\begin{equation*}
\label{eqn:decompositionN}
\psi_N=\psi_N(v,g)=\psi_N(v,0)+\psi_N(0,g):=\psi_N(v)+\tilde{\psi}_N(g)
\end{equation*}

For given Cauchy data $f \in H^{1/2}(\Gamma_O)$ and $g \in H^{-1/2}(\Gamma_O)$ 
we look for $u \in \mathcal{U} = H^{1/2}(\Gamma_I)$ 
satisfying $J(u)=\underset{v \in \mathcal{U}}{\mathrm{inf}} J(v)$ 
where $J$ is the regularized error functional defined by
\begin{equation*}
\label{eqn:kohnvog}
J(v)=\ds \frac{1}{2} \int_{\Omega} \frac{1}{r} 
|| \nabla \psi_D(v,f) - \nabla \psi_N(v,g) ||^2 dx + 
\ds \frac{\varepsilon}{2} \int_\Omega \frac{1}{r} ||\nabla \psi_D(v)||^2 dx
\end{equation*}
In this expression the first term measures a misfit between the Dirichlet solution and the Neumann solution 
whereas the second one is a regularization term in which $\varepsilon$ is a regularization parameter.

Let $u,v \in H^{1/2}(\Gamma_I)$ and define the bilinear symmetric forms
\begin{equation*}
s_D(u,v)= \ds \int_\Omega \ds \frac{1}{r} \nabla \psi_D(u) \nabla \psi_D(v) dx 
\end{equation*}

\begin{equation*}
s_N(u,v)= \ds \int_\Omega \ds \frac{1}{r} \nabla \psi_N(u) \nabla \psi_N(v) dx 
\end{equation*}

and the linear form $l$ 
\begin{equation*}
l(v)=-\ds \int_\Omega \ds \frac{1}{r} (\nabla \tilde{\psi}_D(f) - 
\nabla \tilde{\psi}_N(g))\nabla \psi_D(v)dx
\end{equation*}

as well as the quantity $c$ depending on the Cauchy data only
\begin{equation*}
c=\ds \frac{1}{2} \int_\Omega \ds \frac{1}{r} ||\nabla \tilde{\psi}_D(f) - \nabla \tilde{\psi}_N(g)||^2dx
\end{equation*}

Following the analysis provided in \cite{ACL.B.Faugeras.12.2} 
the functional $J$ can be rewritten as
\begin{equation*}
J(v)=\ds \frac{1}{2}( (1+\varepsilon) s_D(v,v)- s_N(v,v) )-l(v)+c 
\end{equation*}
and its minimum is given by the first order optimality condition or Euler equation which is linear and reads
\begin{equation*}
\label{eqn:opt0}
(J'(u),v)= (1+\varepsilon)s_D(u,v) - s_N(u,v) -l(v)=0\quad \forall v \in \mathcal{U}
\end{equation*}
Solving this equation gives the optimal Dirichlet conditions $u_{opt}$ on $\Gamma_I$ from which one can compute 
$\psi_D(u_{opt},f)$ or $\psi_N(u_{opt},g)$ and then look for the iso-contour defining the plasma boundary.

\section*{\bf III. NUMERICAL IMPLEMENTATION AND RESULTS}
\label{sec:num}

\subsection*{\bf III.A. The toroidal harmonics expansion}

The Legendre functions of degree 1 and half integer order involved in the toroidal harmonic 
expansion are computed 
thanks to the Fortran code provided together with reference \cite{Segura:2000}. 

In the first step of the method one also needs to provide a pole for the toroidal coordinate system and the 
number of harmonics to be used. 

Because the internal toroidal harmonics functions are singular at the pole $F_0$ of the coordinate system, 
the expansion procedure provides meaningful results if the pole lies inside the unknown 
plasma region and not too close from the boundary.
An excellent candidate for the choice of this pole is the current center $(r_c,z_c)$ defined from the moments 
of the plasma current density \cite{Zakharov:1973, Braams:1991}. 
\begin{eqnarray*}
	I_p:=\ds \int_{D_\Gamma} j_p  dx & = & \int_\Gamma \frac{1}{\mu_0} B_s ds \\  
	z_c I_p:=\ds \int_{D_\Gamma} z j_p dx & = & \int_\Gamma \frac{1}{\mu_0}(-r\log r B_n + z B_s) ds \\  
	r^2_c I_p:=\ds \int_{D_\Gamma} r^2 j_p dx & = & \int_\Gamma \frac{1}{\mu_0}(2rz B_n + r^2 B_s) ds 
\end{eqnarray*}
where $\Gamma$ is any closed contour containing the plasma and $D_\Gamma$ the domain it defines.
These quantities can be directly approximated as weighted sums of magnetic measurements 
and then precisely recomputed as integrals on the contour $\Gamma_0$ at every point of which the 
flux $\psi$ and the field $B$ can be evaluated.  

\begin{figure}[h]
\begin{center}
\begin{tabular}{c}
\includegraphics[height=8cm,angle=0]{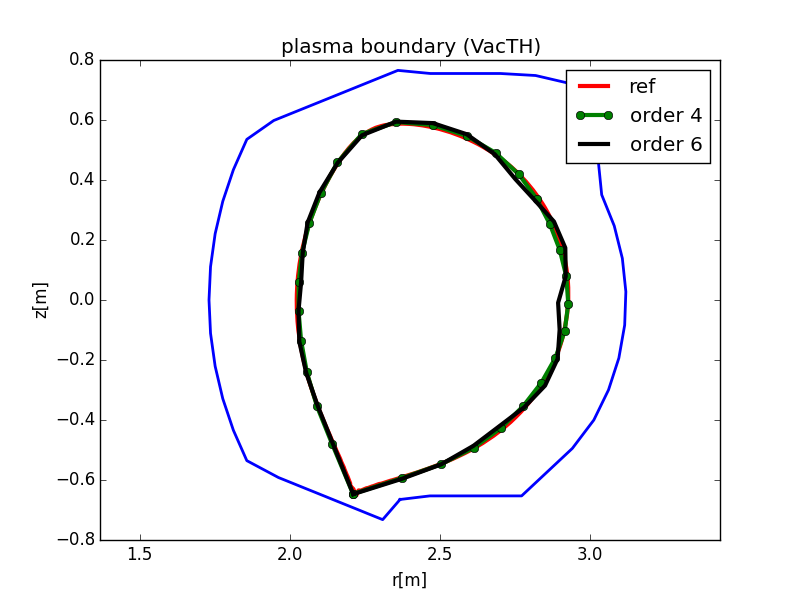} 
\end{tabular}
\end{center}
\caption{\label{fig:west-vacth} 
A WEST equilibrium case. The reference plasma boundary from the equilibrium code CEDRES++ 
and the boundary reconstructed using harmonics of order 4 in VacTH (bullets) 
are almost superimposed. The boundary reconstructed using harmonics of order 6 presents 
some unrealistic oscillations.}
\end{figure}

\begin{figure}[h]
\begin{center}
\begin{tabular}{c}
\includegraphics[height=8cm,angle=0]{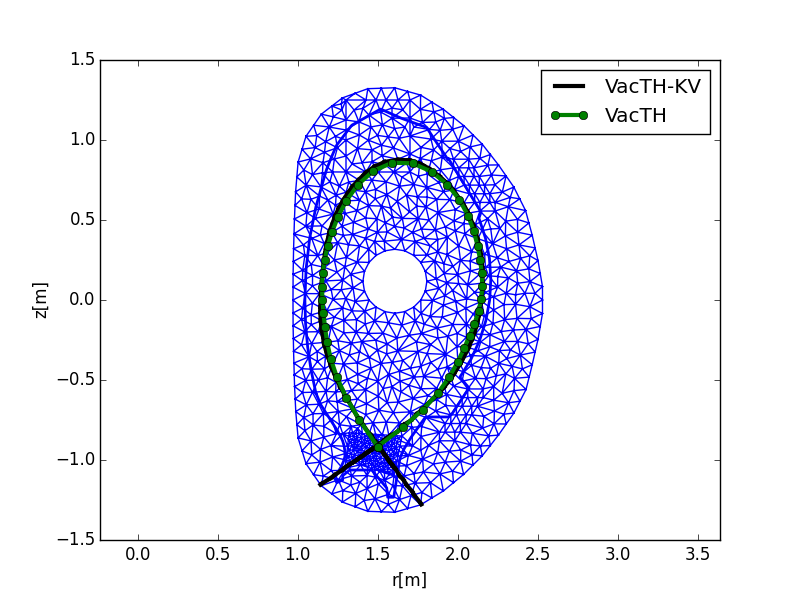} \\
\includegraphics[height=8cm,angle=0]{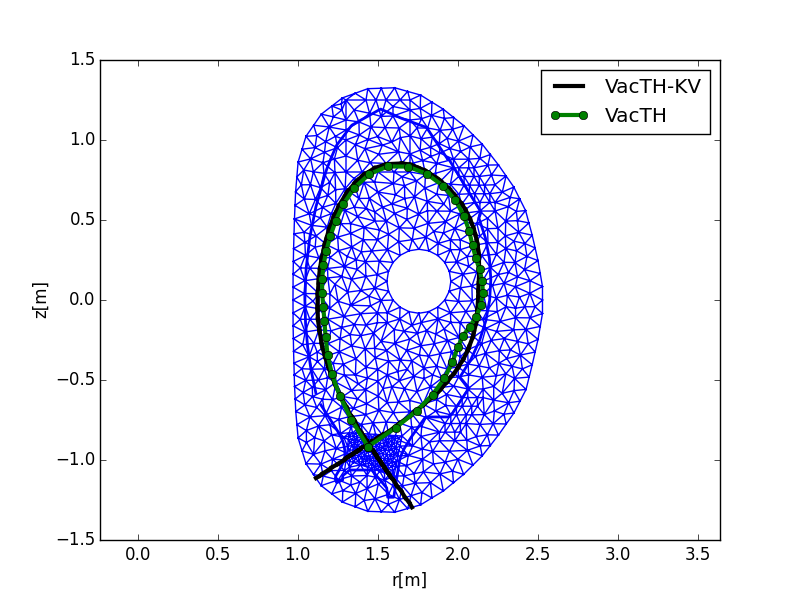}
\end{tabular}
\end{center}
\caption{\label{fig:aug} AUG shot 25374 plasma boundary reconstructions at time $1.1$s (top) and time $2.5$ (bottom).
The plasma boundary computed from toroidal harmonics of order 4 in VacTH (bullets) 
is satisfatory at $1.1$s but presents 
unrealistic oscillations at $2.5$s. 
The plasma boundaries computed with the new proposed method VacTH-KV 
are shown with continuous lines.}
\end{figure}

The number of harmonics to be used need not be large to obtain precise fit 
to the measurements \cite{ACL.B.Faugeras.14.2}. 
Typically the order of the harmonics is chosen to be $4$, $5$ or $6$.
Figure \ref{fig:west-vacth} shows a WEST example case of boundary reconstruction with the code using 
toroidal harmonics only (called VacTH). 
Two reconstructions are performed using harmonics of order 4 (case 1) and 6 (case 2). 
In both cases the root mean square errors for the fit to the measurements is of $2.0\times10^{-3} [T]$ 
for B sensors and 
$1.0 \times 10^{-3} [Wb]$ for flux loops. The Cauchy conditions computed in both cases differ from only small values: 
$\underset{\Gamma_O}{max}|\psi_{case1} - \psi_{case2}|=2.0 \times 10^{-4}$ and 
$\underset{\Gamma_O}{max}|\frac{1}{r}\frac{\partial \psi}{\partial n}_{case1} 
- \frac{1}{r}\frac{\partial \psi}{\partial n}_{case2}|=5.0 \times 10^{-4}$.
However the plasma boundary reconstructed with VacTH differs quite a lot from one case to the other. 
In case 1 it is accurate whereas in case 2 it presents unrealistic oscillations. This shows two things. 
On the one hand the interpolation of the magnetic measurements using toroidal harmonics and the evaluation of 
Cauchy data on the contour $\Gamma_O$ is accurate regardless of the order of the harmonics used. On the other hand 
the reconstructed plasma boundary is quite dependent on the order of the harmonics used and particularly on the 
order of the internal harmonics. For all WEST cases that have been treated up today the choice 
of order 4 harmonics allways gave accurate results. However it can be understood from this example 
that the choice of this order might not be easy or even possible for other Tokamak configurations. 
Figure \ref{fig:aug} shows such an issue. 
We used VacTH to reconstruct the plasma boundary using ASDEX UpGrade experimental data for shot 25374 from 
the ITM-WPCD database \cite{ITM}. The same choice for the order of the toroidal harmonics gives a good reconstruction at 
time $1.1$s but this is not the case at time $2.5$s. This illustrates the need for the optimal control method 
proposed in this paper and implemented in a code called VacTH-KV.

\subsection*{\bf III.B. The optimal control problem}

The resolution of the optimal control problem in VacTH-KV is based on 
a classical $P^1$ finite element method \cite{Ciarlet:1980}. 
Although over choices would have been possible 
the advantage of this numerical method over finite differences or boundary integral methods 
is its simplicity when it comes to adapt the mesh to any internal and external contours, 
and to extract the plasma boundary from the finite dimensional representation of the flux. 

Given the fixed domain $\Omega$ let us consider the family of triangulation ${\tau}_h$ of
$\Omega$, and $V_h$ the finite dimensional subspace of $H^1(\Omega)$
defined by
$$
V_h=\{\psi_h\in H^1(\Omega), \psi_{h|T}\in P^1(T),\,\forall T\in {\tau}_h\}.
$$
Let us also introduce the finite element space on $\Gamma_I$
$$
D_h=\{v_h=\psi_h|_{\Gamma_I},\  \psi_h \in V_h\}.
$$
Consider $(\phi_i)_{i=1,...N}$ a basis of $V_h$ and assume that the first $N_{\Gamma_I}$ 
mesh nodes (and basis functions) correspond to the ones situated on $\Gamma_I$. 
A function $\psi_h \in V_h$ is decomposed as $\psi_h=\sum_{i=1}^N \psi_i \phi_i$ and its trace on 
$\Gamma_I$ as $v_h=\psi_h|_{\Gamma_I}=\sum_{i=1}^{N_{\Gamma_I}} \psi_i \phi_i|_{\Gamma_I}$.

Given boundary conditions $v_h$ on $\Gamma_I$  and $f_h$, $g_h$ on $\Gamma_O$ 
one can compute the approximations 
$\psi_{D,h}(v_h)$,  $\psi_{N,h}(v_h)$, $\tilde{\psi}_{D,h}(f_h)$ and $\tilde{\psi}_{N,h}(g_h)$ 
with the finite element method. These computations involve two different linear systems. One 
for the Dirichlet-Dirichlet problem and one for the Dirichlet-Neumann problem with 
two symmetric positive definite matrices $\mathbf{A}_{DD}$ and $\mathbf{A}_{DN}$.

In order to minimize the discrete regularized error functional, $J_{\varepsilon,h}(u_h)$ 
we have to solve the discrete optimality condition which reads
\begin{equation*}
\label{eqn:optdiscrete}
(1+\varepsilon) s_{D,h}(u_{h},v_h) - s_{N,h}(u_h,v_h) - l(v_h)=0\quad \forall v_h \in D_h
\end{equation*}
which is equivalent to look for the vector $\mathbf{u}$ solution to the linear system
\begin{equation*}
\label{eqn:linearsystem}
\mathbf{Su}=\mathbf{l}
\end{equation*}
where the $N_{\Gamma_I} \times N_{\Gamma_I}$ matrix $\mathbf{S}$ represents 
the bilinear form $s_h=(1+\varepsilon) s_{D,h} -s_{N,h}$ and is defined by
\begin{equation*}
\mathbf{S}_{ij}= s_h(\phi_i,\phi_j)
\end{equation*}
and $\mathbf{l}$ is the vector $(l_h(\phi_i))_{i=1,...N_{\Gamma_I}}$.

The matrices which constitues $\mathbf{S}$ are evaluated by
\begin{equation*}
s_{D,h}(\phi_i,\phi_j)=\ds \int_\Omega \frac{1}{r}\nabla \psi_{D,h}(\phi_i) \nabla \mathcal{R} (\phi_j) dx
\end{equation*} 
and 
\begin{equation*}
s_{N,h}(\phi_i,\phi_j)=\ds \int_\Omega \frac{1}{r}\nabla \psi_{N,h}(\phi_i) \nabla \mathcal{R} (\phi_j) dx
\end{equation*} 
where $\mathcal{R}(\phi_j)$ is the trivial extension which coincides with $\phi_j$ on $\Gamma_I$ 
and vanishes elsewhere.

In the same way the right hand side $\mathbf{l}$ is evaluated by
\begin{equation*}
l_h(\phi_i)=-\ds \int_\Omega \frac{1}{r} ( \nabla \tilde{\psi}_{D,h}(f_h)- \nabla \tilde{\psi}_{N,h}(g_h)) 
\nabla \mathcal{R} (\phi_i) dx
\end{equation*}

It should be noticed here that matrices $\mathbf{A}_{DD}$, $\mathbf{A}_{DN}$ and $\mathbf{S}$ 
depend on the geometry of the problem only. 
Therefore an internal contour $\Gamma_I$ being given 
they can be computed once as well as their Cholesky decomposition 
and be used for the resolution of successive problems with varying input data. 
Only the righthand side $\mathbf{l}$ has to be recomputed. 

\begin{figure}[h]
\begin{center}
\begin{tabular}{c}
\includegraphics[height=8cm,angle=0]{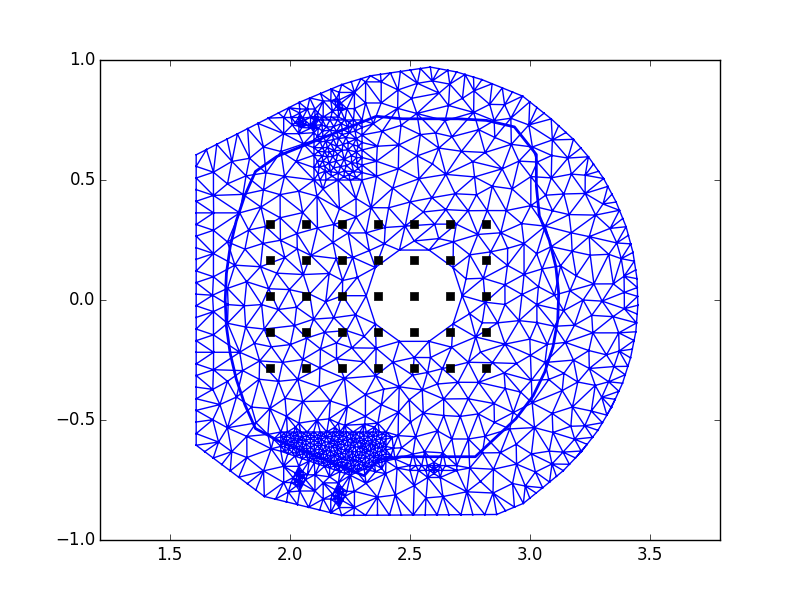}\\
\includegraphics[height=8cm,angle=0]{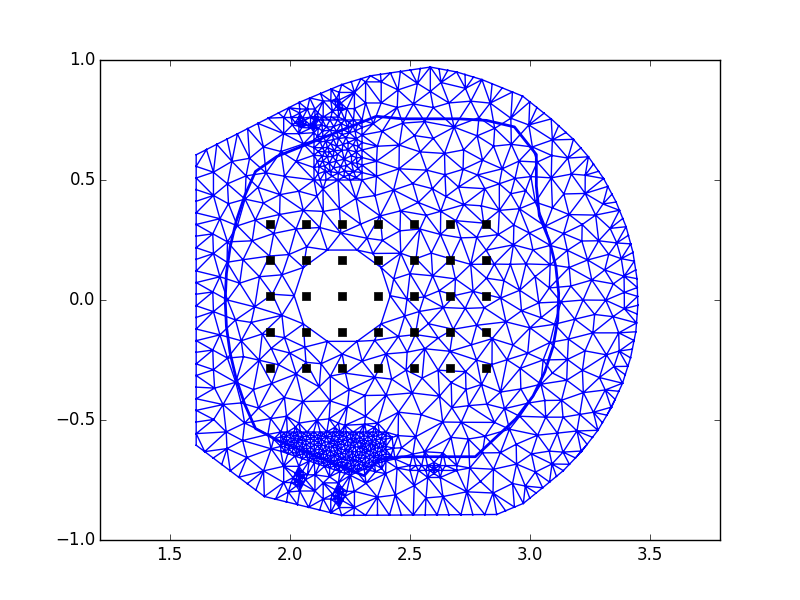}
\end{tabular}
\end{center}
\caption{\label{fig:westholes} Examples of precomputed meshes for WEST. 
The limiter is represented. The mesh is refined in lower and upper X-point regions.
Each small square inside the limiter region is the center of a circle defining an internal 
contour $\Gamma_I$. }
\end{figure}

An important point however is that in the same way that the pole $F_0$ of the toroidal 
coordinate system moves following the plasma current center, the inner contour $\Gamma_I$ might also need 
to be displaced. We have adopted the following strategy consisting in defining a number of points in 
the vacuum vessel. Each of them is the center of a circle defining a contour $\Gamma_I$. A different mesh is precomputed 
for each of these contours as well as the associated finite elements matrices (see Figure \ref{fig:westholes} 
for an example). 
Each time a new set of magnetic measurements is provided the plasma current center is computed in the first step 
of the algorithm involving the toroidal harmonics expansion. Then one looks for the internal contour whose center 
is the closest to this point and use the associated and precomputed mesh and matrices in the second step 
in order to minimize the Kohn-Vogelius functional. This enables fast computations. The computing time is of the order 
of $1ms$ for one boundary reconstruction in the WEST configuration on a Laptop with two quadcore 
processors of frequency 2.40 Ghz.

\begin{figure}[h]
\begin{center}
\begin{tabular}{c}
\includegraphics[height=8cm,angle=0]{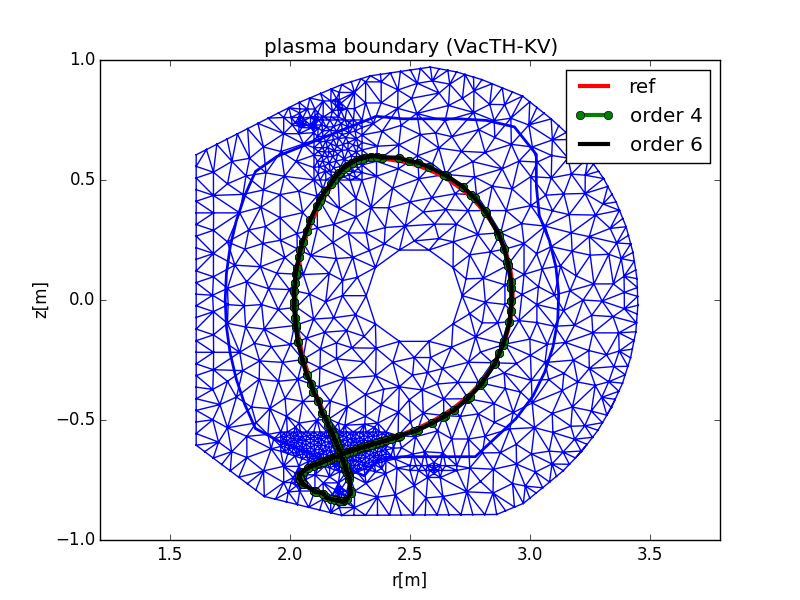}
\end{tabular}
\end{center}
\caption{\label{fig:west-vacthkv}
Same WEST test case as in Figure \ref{fig:west-vacth}. The boundaries reconstructed with the new 
optimal control proposed method, VacTH-KV, using toroidal harmonics of order either 4 or 6 are almost superimposed 
with the reference boundary from CEDRES++}
\end{figure}

\begin{figure}[h]
\begin{center}
\begin{tabular}{c}
\includegraphics[height=8cm,angle=0]{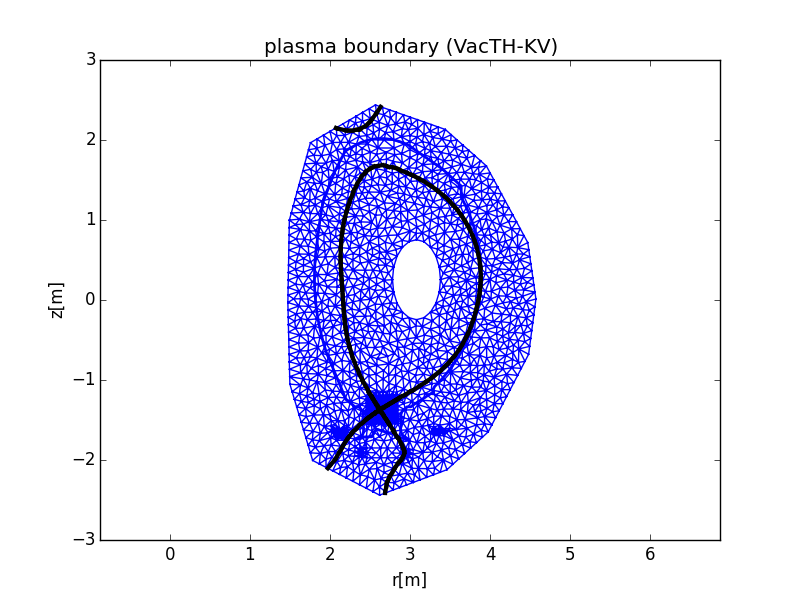}
\end{tabular}
\end{center}
\caption{\label{fig:jet-vacthkv}
VacTH-KV boundary reconstruction for a JET test case example (shot 74221 at 46s). The method using toroidal harmonics only, VacTH,
could not reconstruct a proper boundary.}
\end{figure}

Figure \ref{fig:west-vacthkv} shows the reconstructed plasma boundaries for 
the same WEST example case as in figure \ref{fig:west-vacth} but this time with VacTH-KV. 
Using toroidal harmonics of order either 4 or 6 for the evaluation of Cauchy conditions 
on $\Gamma_O$ leads to the same reconstruction almost superimposed with the reference boundary from 
the equilibrium code CEDRES++. 
In addition to the unrealistic results given by VacTH for the two AUG cases with real experimental 
measurements, figure \ref{fig:aug} also shows the excellent results provided by VacTH-KV on these 
same cases. As a last example illustrating the non tokamak dependance of the method 
figure \ref{fig:jet-vacthkv} shows the plasma boundary reconstructed by VacTH-KV 
on a JET case (shot 74221 at time 46s). In this case VacTH alone was not able to reconstruct a proper 
plasma boundary.

In all our numerical experiments the regularization parameter $\varepsilon$ is given the value 
of $5 \times 10^{-4}$ as suggested by the L-curve method analysis \cite{Hansen:1998} conducted in 
\cite{ACL.B.Faugeras.12.2}.

Finally we would like to point out the fact that the optimal control method 
presented in section II.B. 
is completely unchanged if poloidal field coils or passive structures with measured currents are included 
in the domain $\Omega$. In fact this is the case for the WEST, AUG and JET numerical results presented here.

\section*{\bf IV. SUMMARY AND CONCLUSION}

In this paper a method for plasma boundary reconstruction is proposed. 
It is decomposed in two steps. A toroidal harmonics expansion is first used to interpolate 
discrete magnetic measurements on an external contour. This interpolation using functions which are 
exact analytic solutions of the equation satisfied by the flux is accurate. On the other hand 
we have shown that the 
extrapolation of the flux towards the plasma boundary using toroidal harmonics only can in some cases be 
inaccurate. This is the reason why we introduce a second step in which a Cauchy problem on a fixed domain is 
solved by an optimal control method. These two steps together provide an accurate and robust plasma boundary 
reconstruction method. It is generic and can be used for any Tokamak having ferromagnetic structures or not. 
An implementation using a finite elements discretization is proposed. The precomputation of several meshes and 
all associated linear algebra quantities allows fast computation. A code called VacTH-KV 
has been developped and is available on the ITM-WPCD platform. 
This paper does not deal with the computation of error bars on the plasma boundary reconstruction results. 
This is of course of interest as for any ill-posed inverse problem and might be conducted in the future for 
example using the epsilon-nets technique \cite{Zaitsev:2014}

\section*{\bf ACKNOWLEDGMENTS}
This work, supported by the European Communities under the contract of Association
between EURATOM-CEA was carried out within the framework of the Task Force on
Integrated Tokamak Modelling of the European Fusion Development Agreement (ITM-WPCD). The
views and opinions expressed herein do not necessarily reflect those of the European
Commission.

The author would like to thank his colleagues at the university of Nice, 
Jacques Blum, Cedric Boulbe and Holger Heumann for many useful discussions. 
The input magnetic measurements for WEST were provided by colleagues from 
the IRFM-CEA at Cadarache, Eric Nardon and Remy Nouailletas, 
and for AUG and JET by the ITM-WPCD project.

Finally the author would like to thank the two anonymous reviewers for their constructive criticism from 
which the paper has benefitted significantly.


\clearpage
\bibliographystyle{ans}
\bibliography{biblio}

\newpage

\newpage

\newpage

\newpage

\newpage

\newpage

\end{document}